\documentclass[11 pt]{amsart}
\usepackage{amsmath}
\usepackage{amssymb}
\usepackage{amsthm}
\usepackage{dsfont}
\usepackage{fancyhdr}
\usepackage[all]{xy}
\usepackage{hyperref}
\usepackage{paralist}
\usepackage{tikz}
\usepackage{verbatim}
\usepackage[hmargin=1.2in,vmargin=1.2in]{geometry}

\newtheorem{theorem}{Theorem}[section]
\newtheorem{lemma}[theorem]{Lemma}
\theoremstyle{proposition}

\theoremstyle{corollary}

\theoremstyle{definition}

\numberwithin{equation}{section}

\newcommand{\Z}{\mathbb Z}

\newcommand{\F}{\mathbb F}

\newcommand{\f}{\mathfrak f}

\newcommand{\g}{\mathfrak g}

\begin{document}
\title[]{On the last fall degree of Weil descent polynomial systems}
\author[]{Ming-Deh A. Huang (USC, mdhuang@usc.edu)}
\address{Computer Science Department,University of Southern California, U.S.A.}
\email{mdhuang@usc.edu}

\urladdr{}
\date{\today}

\maketitle

\begin{abstract}
Given a polynomial system $\mathcal{F}$ over a finite field $k$ which is not necessarily of dimension zero, we consider
the Weil descent $\mathcal{F}'$ of $\mathcal{F}$ over a subfield $k'$.
We prove a theorem which relates the last fall degrees of $\mathcal{F}_1$ and $\mathcal{F}'_1$, where
the zero set of $\mathcal{F}_1$ corresponds bijectively to
the set of $k$-rational points of $\mathcal{F}$, and the zero set of $\mathcal{F}'_1$ is the set of $k'$-rational points of the Weil descent $\mathcal{F}'$.
As an application we derive upper bounds on the last fall degree of $\mathcal{F}'_1$ in the case where $\mathcal{F}$ is a set of linearized polynomials.

\end{abstract}
\section{Introduction}
Let $k$ be a field and let $\mathcal{F} \subset R=k[X_0,\ldots,X_{m-1}]$ be a finite subset which generates an ideal. Let $R_{\leq i}$ be the set of polynomials in $R$ of degree at most $i$.

For $i \in \Z_{\geq 0}$, we let $V_{\mathcal{F},i}$ be the smallest $k$-vector space of $R_{\leq i}$ such that
\begin{enumerate}
\item
$\{f \in \mathcal{F}: \deg(f) \leq i \} \subseteq V_{\mathcal{F},i}$;
\item
if $g \in V_{\mathcal{F},i}$ and if $h \in R$ with $\deg(hg) \leq i$, then $hg \in V_{\mathcal{F},i}$.
\end{enumerate}

We write $f\equiv_i g \pmod{\mathcal{F}}$, for $f,g\in R$, if $f-g\in V_{\mathcal{F},i}$.

The \emph{last fall degree} as defined in \cite{HKY} (see also \cite{HKYY}) is the largest $d$ such that $V_{\mathcal{F},d} \cap R_{\leq d-1} \neq V_{\mathcal{F},d-1}$. We denote the last fall degree of $\mathcal{F}$ by $d_{\mathcal{F}}$.

As shown in \cite{HKY,HKYY} the last fall degree  is intrinsic to a polynomial
system, independent of the choice of a monomial order, always bounded by the degree of regularity,  and  invariant under linear change of variables and linear change of equations.
In \cite{HKY, HKYY} complexity bounds on solving zero dimensional polynomial systems were proven based on
the last fall degree.  It was shown in \cite{HKY} that the polynomial systems arising from
the Hidden Field Equations (HFE) public key crypto-system \cite{BET,DIN} have bounded last fall degree  if
the degree of the defining polynomial and the cardinality of the base
field are fixed (the bound was improved in \cite{GMP}), and
it follows that
the HFE polynomials systems can be solved unconditionally in polynomial time.

For $\mathcal{F} \subset R=k[X_0,\ldots,X_{m-1}]$, let $Z_k (\mathcal{F})$ denote the set of solutions of $\mathcal{F}$ over $k$; let $Z (\mathcal{F})$ denote the set of solutions of $\mathcal{F}$ over $\overline{k}$, where $\overline{k}$ is an algebraic closure of $k$.
If $\mathcal{F}$ is zero-dimensional then determining $Z (\mathcal{F})$ reduces to computing $V_{\mathcal{F},\max(d_{\mathcal{F}},e)}$ where $e$ is the cardinality of $Z (\mathcal{F})$ \cite{HKYY}.

Suppose $k$ is a finite field of cardinality $q^n$ with subfield $k'$ of cardinality $q$.
The Weil descent system of $\mathcal{F}$ to $k'$ is a polynomial system obtained when one expresses all equation with the help of a basis of $k'/k$.   Let $\alpha_0,\ldots,\alpha_{n-1}$ be a basis of $k/k'$. For $f \in \mathcal{F}$ and $j=0,\ldots,n-1$, we define $f_{j} \in k'[X_{ij}, i=0,\ldots,m-1, j=0,\ldots,n-1]$ by
\begin{eqnarray*}
f\left(\sum_{j=0}^{n-1} \alpha_j X_{0j}, \ldots, \sum_{j=0}^{n-1} \alpha_j X_{m-1\ j} \right) = \sum_{j=0}^{n-1} f_{j} \alpha_j.
\end{eqnarray*}
We note that $\deg f_j \le \deg f$.
The system
\begin{eqnarray*}
\mathcal{F}'=\{f_j:\ f \in \mathcal{F}, j=0,\ldots,n-1\}
\end{eqnarray*}
is called the \emph{Weil descent system} of $\mathcal{F}$ with respect to $\alpha_0,\ldots,\alpha_{n-1}$.

There is a bijection between $Z_k (\mathcal{F})$ and $Z_{k'} (\mathcal{F}')=Z({\mathcal{F}}'_1)$, where ${\mathcal{F}}'_1$ is $\mathcal{F}'$ together with the field  equations of $k'$, that is,
$$\mathcal{F}'_1 =  \mathcal{F}'\cup\{  X_{ij}^q-X_{ij}, i=0,\ldots,m-1, j=0,\ldots,n-1 \}.$$

The HFE polynomial system is constructed by forming the Weil descent of some $\mathcal{F}$ consisting of a single univariate polynomial, followed by linear change of variables and linear change of equations \cite{HKY}. Multivariate-HFE systems can  be constructed similarly except $\mathcal{F}$ is replaced by a finite set of multivariate polynomials of dimension zero.  In \cite{HKYY} upper bounds on the last fall degree degree of ${\mathcal{F}}'_1$ were proven in terms of $q$, $m$, the last fall degree of $\mathcal{F}$, the degree of $\mathcal{F}$ and the number of solutions of $\mathcal{F}$, but not on $n$.  The result implies that multi-HFE cryptosystems giving rise to multi-HFE polynomial systems as described above are vulnerable to attack as well.

In this paper we consider the situation where $\mathcal{F}$ is not necessarily zero-dimensional.

Let
\begin{eqnarray*}
\mathcal{F}_1 =   \mathcal{F}\cup\{ X_i^q - Y_{i1},\ldots, Y_{i \ n-2}^q - Y_{i \ n-1}, Y_{i \ n-1}^q - X_i : i=0,\ldots, m-1\}\\
\subset k[X_i, Y_{ij}: i=0,\ldots, m-1; j=1,\ldots, n-1] .
\end{eqnarray*}

We observe that
$Z_k (\mathcal{F})$ can easily be identified with  $Z(\mathcal{F}_1)$.
So there is a bijection between $Z(\mathcal{F}_1)$ and $Z({\mathcal{F}}'_1)$.
Note also that the ideals generated by $\mathcal{F}_1$ and $\mathcal{F}'_1$ are radical ideals.

The following theorem relates the last fall degrees of $\mathcal{F}_1$ and $\mathcal{F}'_1$.
\begin{theorem}
\label{weil}
$\max( d_{\mathcal{F}_1}, q \deg \mathcal{F} ) = \max ( d_{\mathcal{F}'_1}, q \deg \mathcal{F} )$.
\end{theorem}

Theorem~\ref{weil} is closely related to Proposition 2 of \cite{HKY} and Proposition 4.1 of \cite{HKYY}.  In comparison, the bound established in Theorem~\ref{weil} is a bit weaker.  However the set $\mathcal{F}_1$ stated in the theorem directly contains $\mathcal{F}$ as a subset.  This makes it easier to apply the theorem both conceptually and technically. When $\mathcal{F}$ consists of a univariate polynomial or more generally when $Z(\mathcal{F})$ is finite, it is not hard to bound $d_{\mathcal{F}_1}$.  From this an easier and more conceptual proof of the theorems in \cite{HKY,HKYY} can be constructed based on Theorem~\ref{weil}.  However in this paper we will focus on applying the theorem to the situation where $\mathcal{F}$ is not zero-dimensional, especially when $\mathcal{F}$ consists of linearized polynomials.

\subsection{Proof of Theorem~\ref{weil}}
For non-negative integers $i$, let $\sigma_i$ denote the automorphism of $\overline{k}$ over $k'$ such that $\sigma_i (x) = x^{q^i}$ for $x\in \overline{k}$.  For every multivariate polynomial $h$ with coefficients from $\overline{k}$, let $h^{\sigma_i}$ denote the polynomial obtained from $h$ by acting on each coefficient of $h$ by $\sigma_i$.

Let $\Gamma$ be the $n$ by $n$ matrix with rows and columns indexed by $0,\ldots, n-1$, so that $\alpha_j^{\sigma_i}$ is the $(i,j)$-th entry of $\Gamma$ for $i,j=0,\ldots,n-1$.

Let
$$g_f=f\left(\sum_{j=0}^{n-1} \alpha_j X_{0j}, \ldots, \sum_{j=0}^{n-1} \alpha_j X_{m-1\ j} \right) = \sum_{j=0}^{n-1} f_{j} \alpha_j$$
where $f_{j} \in k'[X_{ij}, i=0,\ldots,m-1, j=0,\ldots,n-1]$.

Let
$\hat{f}=\left( \begin{array}{c} f_0\\.\\.\\.\\f_{n-1}\end{array}
\right)$ be the column vector with
$f_i$ as the $i$-th entry for $i=0,\ldots,n-1$.

Then $g_f^{\sigma_i}=\sum_{j=0}^{n-1} f_{j} \alpha_j^{\sigma_i}$, and
$\Gamma \hat{f} = \left( \begin{array}{c} g_f^{\sigma_0}\\.\\.\\.\\g_f^{\sigma_{n-1}}\end{array}
\right)$.
Let $\mathcal{G}=\{ g_f^{\sigma_0},\ldots,g_f^{\sigma_{n-1}}: f\in\mathcal{F}\}$ and
$\mathcal{G}_1 = \mathcal{G}\cup \{  X_{ij}^q-X_{ij}, i=0,\ldots,m-1, j=0,\ldots,n-1 \}$.
Since $\Gamma$ is invertible, it follows from Proposition 2.6 (part v) of \cite{HKYY} that $d_{\mathcal{G}} = d_{\mathcal{F'}}$, and
$d_{\mathcal{G}_1} = d_{\mathcal{F'}_1}$.

Let $Z_{ij}$, $i=0,\ldots,m-1$ and $j=0,\ldots,n-1$, be defined by the following change of coordinates:
$$\left( \begin{array}{c} X_i\\Y_{i1}\\.\\.\\Y_{i \ n-1}\end{array}
\right)
= \Gamma
\left( \begin{array}{c} Z_{i0}\\Z_{i1}\\.\\.\\Z_{i \ n-1}\end{array}
\right).$$
Under the change of coordinates, $\mathcal{F}_1$ becomes $\mathcal{G}_2$ where
\begin{eqnarray*}
\mathcal{G}_2 & = & \{ f (\sum_{j=0}^{n-1} \alpha_j Z_{0j},\ldots,\sum_{j=0}^{n-1} \alpha_j Z_{m-1\ j}): f\in\mathcal{F}\}\cup\{ Z_{ij}^q- Z_{ij}:i,j=0,\ldots,n-1\}\\
& = & \{ g_f(Z_{01},\ldots,Z_{m,n-1}): f\in\mathcal{F}\}\cup\{ Z_{ij}^q- Z_{ij}:i,j=0,\ldots,n-1\},
\end{eqnarray*}
which we identify as a subset of $\mathcal{G}_1$.
Since $g_f^q \equiv g_f^{\sigma} \mod I$ where $I$ is the deal generated by $X_{ij}^q - X_{ij}$, $i=0,\ldots, m-1$, $j=0,\ldots, n-1$,  we see that $g_f^{\sigma}\in V_{\mathcal{G}_2, qd}$ where $d=\deg\mathcal{F} \ge \deg g_f$.  It follows inductively that $g_f^{\sigma_i}\in V_{\mathcal{G}_2, qd}$ for all $i$, hence $\mathcal{G}_1 \subset V_{\mathcal{G}_2, qd}$.  Hence $V_{\mathcal{G}_1, i}=V_{\mathcal{G}_2, i}$ for $i\ge qd$.  Therefore
$\max(d_{\mathcal{G}_1},qd)=\max(d_{\mathcal{G}_2},qd)$.  Since
$d_{\mathcal{G}_1}=d_{\mathcal{F}'_1}$ and $d_{\mathcal{F}_1}=d_{\mathcal{G}_2}$, we conclude that
$\max(d_{\mathcal{F}_1},qd)=\max(d_{\mathcal{F'}_1},qd)$.  Theorem~\ref{weil} follows.

\section{Systems of linearized polynomials}
A $k'$-linearized polynomial in $R=k[x_0\ldots x_{m-1}]$ is an element of the $k$-submodule of $R$ generated by
$x_i^{q^j}$ where $q=|k'|$, $i=0\ldots m-1$ and $j\ge 0$.   As before let $n=[k:k']$.  Let $Q = \{ x_i^{q^n}-x_i : i=0\ldots m-1\}$.

For $g=\sum_{i=0}^{d} a_i x^i\in k[x]$,  let $L(g)=\sum_{i=0}^d a_i x^{q^i}$.  More generally we consider the $k$-linear map from $\oplus_{i=0}^{m-1} k[x_i]$ onto the $k$-module of $k'$-linearized polynomials such that $L(x_i^j)=x_i^{q^j}$ for $i=0,\ldots,m-1$ and $j\ge 0$.

Let $S=k [x_{ij}: i=0\ldots m-1, j=0,\ldots,n-1]$.  We also write $S=k[\hat{x}_i: i=0,\ldots,m-1]$. where
$\hat{x}_i = x_{i0},\ldots,x_{i\ n-1}$.  Let $S_1\subset S$ be the $k$-module of linear forms over $x_{ij}$, $i=0\ldots m-1, j=0,\ldots,n-1$.

For $g\in k[x]$ and $h\in k[x_i]$, let $g\circ h \in k[x_i]$ be defined as $(g\circ h) (x_i) = g(h(x_i))$. For $f=\sum_{i=0}^{m-1} f_i \in\bigoplus_{i=0}^{m-1} k[x_i]$ with $f_i\in k[ x_i]$, let $g \circ f\in \bigoplus_{i=0}^{m-1} k[x_i]$ be defined as $g\circ f = \sum_{i=0}^{m-1} g\circ f_i$.  Hence $(g\circ f) (x_0,\ldots,x_{m-1}) =  \sum_{i=0}^{m-1} g(f_i (x_i))$.  Note also that $L(g)\circ L(f)$ is a $k'$-linearized polynomial in $k[x_0,\ldots,x_{m-1}]$.  

Similarly for $g\in k[x]$ and $f=\sum_{i,j} f_{ij} \in\bigoplus_{i=0}^{m-1}\bigoplus_{j=0}^{n-1} k[x_{ij}]$
with $f_{ij}\in k[ x_{ij}]$, let $g \circ f\in \bigoplus_{i=0}^{m-1}\bigoplus_{j=0}^{n-1} k[x_{ij}]$ be defined as $g\circ f = \sum_{i,j} g\circ f_{ij}$.

Consider the map $\ell$ from $\bigoplus_{i=0}^{m-1}\bigoplus_{j=0}^{n-1} k[x_i^j]$ to  $S_1$
such that $\ell (x_i^j) = x_{ij}$.

Let $\bar{Q}=\{x_{ij}^q-x_{i\ j+1}: i=0\ldots m-1, j=0\ldots n-1\}$ where $j+1$ is taken $\mod n$.
Consider the $k$-algebra isomorphism from $R/\langle Q \rangle$ to $S/\langle \bar{Q} \rangle$ sending $x_i^{q^j}$ to $x_{ij}$ for
$i=0,\ldots,m-1$, $j=0,\ldots,n-1$.  (Note that $x_i^{q^{j+1}} = (x_i^{q^j})^q$ maps to $x_{ij}^q$ and  $x_{ij}^q \equiv x_{i\ j+1} \mod \bar{Q}$.)

For $f\in R$ where the degree of $f$ in $x_i$ is less than $q^n$ for all $i$, let $\bar{f}\in S$ denote the image of $f$ in  $S/\langle \bar{Q} \rangle$  under the isomorphism.
  We note that
elements of $S_1$ are all distinct mod $\bar{Q}$.   Let $f=\sum_{i=0}^{m-1} f_i (x_i)$ with $\deg f_i < n$ for all $i$.
Let $f_i = \sum_{j=0}^{n-1} a_{ij} x_i^j$.  Then $L(f)\in R$ corresponds to $\ell(f)$.  If we identify with $x_i\in R$ with $x_{i0}\in S$ for $i=0,\ldots,m-1$.  Then $L(f)\equiv_d \ell(f)\pmod{\bar{Q}}$ where $d=\deg L(f)$.

Suppose $\mathcal{F}$ is a finite set of $k'$-linearized polynomials of maximum degree $d=q^c$ for some $c > 0$.
We may identify $x_i\in R$ with $x_{i0}\in S$ and consider $\mathcal{F}\subset S$.  Let $\mathcal{F}'$ be the Weil descent system of $\mathcal(F)$ with respect to a $k/k'$ basis.  We are interested in the last fall degree of $\mathcal{F}'_1=\mathcal{F}'\cup\{ x_{ij}^q-x_{ij}: i=0,\ldots,m-1, j=0,\ldots,n-1\}$.
Let $\mathcal{G} =\mathcal{F}\cup \bar{Q}\subset S$. By Theorem\ref{weil}, $\max( d_{\mathcal{G}}, q \deg \mathcal{F} ) = \max ( d_{\mathcal{F}'_1}, q \deg \mathcal{F} )$.

Recall that $Z(\mathcal{F}'_1)=Z_{k'} (\mathcal{F}')$, which corresponds to $Z_{k} (\mathcal{F})=Z(\mathcal{F}\cup Q )$, the set of $k$-rational points of $Z(\mathcal{F})$.  In what follows we consider a more general situation where instead of $Z_{k} (\mathcal{F})$ we are interested in $Z_{W} (\mathcal{F})= Z(\mathcal{F})\cap W^{m}$ where $W$ is a $\tau$-invariant subspace of $k$ and $\tau$ is the Frobenius map over $k'$: $x \to x^q$ for all $x\in k$.   Note that every $\tau$-invariant subspace $W$ of $k$ is of of the form $Z(L(\f_W))$ where $\f_W$ divides $x^n-1$.  In fact $W$ is the kernel of $\f_W(\tau)$, and $f_W$  is the characteristic polynomial of $\tau$ as a linear map on $W$.  In particular $\f_W=x^n-1$ corresponds to $W=k$ and $\f_W=x-1$ corresponds to $W=k'$.  Suppose $d_W = \deg \f_W$.

In this more general situation we let
$S=k [x_{ij}: i=0\ldots m-1, j=0,\ldots,d_W -1]$.  We also write $S=k[\hat{x}_i: i=0,\ldots,m-1]$. where
$\hat{x}_i = x_{i0},\ldots,x_{i\ d_W-1}$.
Let $f=\sum_{i=0}^{m-1} f_i (x_i)$ with $\deg f_i < d_W$ for all $i$.
Suppose $f_i = \sum_{j=0}^{d_W-1} a_{ij} x_i^j$.  Then $L(f)=\sum_{i=1}^{m-1} L(f_i)$ where $L(f_i)=\sum_{j=0}^{d_W -1} a_{ij} x_i^{q^j}$ and  $\ell(f)=\sum_{i=0}^{m-1} \ell(f_i)$ where $\ell (f_i) = \sum_{j=0}^{d_W-1} a_{ij} x_{ij}$.

Below we fix $W$ and let $n'=d_W$.  Write $\f_W (x) = x^{n'} - \g_W (x)$ with $\deg \g_W < n'$.
Let $Q=\{x_i^{q^{n'}} - L(\g_W(x_i)): i=0,\ldots,m-1\}$, and correspondingly we let
$\bar{Q}=\{x^q_{i \ n'-1} -  \ell (\g_W (x_i)), x_{ij}^q - x_{i \ j+1}: i=0,\ldots,m-1, j=0,\ldots,n' -2\}$.
Then we
have an isomorphism from $R/\langle Q \rangle$ to $S/\langle \bar{Q} \rangle$ sending $x_i^{q^j}$ to $x_{ij}$ for
$i=0,\ldots,m-1$, $j=0,\ldots, n' -1$. Let $S_1\subset S$ be the $k$-module of linear forms over $x_{ij}$, $i=0\ldots m-1, j=0,\ldots,n' -1$.  We note that
elements of $S_1$ are all distinct mod $\bar{Q}$.

For $f\in R$, we have $f\equiv_d f_1 \pmod{Q}$ where $d=\deg f$ and the degree of $f_1$ in $x_i$ is less than $q^{d_W}$ for all $i$.  Let $\bar{f}\in S$ denote the image of $f_1$ in  $S/\langle \bar{Q} \rangle$  under the isomorphism.
Let $f=\sum_{i=0}^{m-1} f_i (x_i)$ with $\deg f_i < n'$ for all $i$.
Then $\overline{L(f)} = \ell(f)$.
If we identify with $x_i\in R$ with $x_{i0}\in S$ for $i=0,\ldots,m-1$.  Then $L(f)\equiv_d \ell(f)\pmod{\bar{Q}}$ where $d=\deg L(f)$.

\begin{lemma}
\label{fiq}
Suppose $f=\sum_{i,j} a_{ij} x_{ij}\in S_1$ with $a_{ij}\in k$. Then with respect to $\bar{Q}$, $f_i^q \equiv_q f_{i+1}\pmod{\bar{Q}}$, for $i\ge 0$ where $f=f_0$ and $f_i\in S_1$ for $i\ge 0$.
\end{lemma}
\ \\{\bf Proof} For $r\ge 0$, we have inductively $f_r = \sum b_{ij} x_{ij}\in S_1$.  Now
$f_r^q = \sum_{ij} b_{ij}^{q} x_{ij}^q$, and since for all $i$, $x_{ij}^q \equiv_q x_{i\ j+1} \pmod{\bar{Q}}$ for $j=0,\ldots, n'-2$, and $x^q_{i \ n'-1}\equiv_q  \ell (\g (x_i))$, it follows that $f_r^q \equiv_q f_{r+1} \pmod{\bar{Q}}$ with $f_{r+1}\in S_1$.  $\Box$

\begin{lemma}
\label{L}
Let $\mathcal{H}$ be a finite set of $S$ and suppose $\bar{Q}\subset \mathcal{H}$.
Suppose $f\in S_1$ and $f\equiv_i 0 \pmod{\mathcal{H}}$ for some $i > 0$.  Let $g\in k[x]$.  Then $L(g)\circ f -  f'\in \langle \bar{Q} \rangle$ for some $f'\in S_1$, where $\langle \bar{Q} \rangle$ denotes the ideal generated by $\bar{Q}$, and $f'\equiv_r 0 \pmod{\mathcal{H}}$ where $r=\max(i,q)$.
\end{lemma}
\ \\{\bf Proof}  The lemma follows by applying Lemma~\ref{fiq} inductively.  More specifically assume inductively $L(x^i)\circ f\equiv_q f'_i \pmod{\bar{Q}}$ with $f'_i\in S_1$, then
$L(x^{i+1})\circ f \equiv_q {f'}^q_i \equiv_q f'_{i+1} \pmod{\bar{Q}}$ for some $f'_{i+1}\in S_1$.  From this the lemma easily follows. $\Box$

For $f\in\mathcal{F}$, $ f\equiv_d \bar{f}\pmod{\bar{Q}}$ with $\bar{f}\in S_1$.  Let $\bar{F}$ consist of all such $\bar{f}\in S_1$ with $f\in\mathcal{F}$.  Let $\bar{\mathcal{G}}=\bar{\mathcal{F}}\cup\bar{Q}$.  Then $\bar{\mathcal{G}}\subset V_{\mathcal{G}, d}$ and $\bar{\mathcal{G}}\subset V_{\bar{\mathcal{G}}, q}$.

Let $S_{1i}=S_1\cap k[\hat{x}_{j}: j=i,\ldots,m-1]$, that is , the submodule containing all $k$-linear forms
in $x_{ij}$, $i=i,\ldots,m-1$, $j=0,\ldots,n' -1$.  In particular $S_1 = S_{10}$.
Let $\bar{Q}_r =\bar{Q}\cap k[ \hat{x}_{i}: i=r,\ldots, m-1]$ for $r=1,\ldots,m-1$.
\begin{lemma}
\label{xn-1}
Consider a $k'$-linearized polynomial of the form $L(f)$ with $f=\sum_{i=0}^{m-1} f_i$ and $f_i \in k[x_i]$ of degree less than $n'$, for $i=0,\ldots,m-1$.
Suppose $\ell(f)\in V_{\bar{\mathcal{G}},q}$ and the GCD of $f_0$ and $\f_W$ is 1.
Then $x_{00} - \ell_0 \in  V_{\bar{\mathcal{G}},q}$ for some linear form $\ell_0 \in S_{11}$.  Moreover for $i=1,\ldots,n' -1$, $x_{0i} - \ell_i \in  V_{\bar{\mathcal{G}},q}$ for some linear form $\ell_i \in S_{11}$, and $\ell_{i-1}^q \equiv \ell_i\pmod{\bar{Q}_{1}}$.
\end{lemma}

\ \\{\bf Proof} Since the GCD of $f_0$ and $\f_W$ is 1,  $A(x) f_0 (x)+ B(x) \f_W (x) =1$ for some $A(x), B(x)\in k[x]$ where $\deg A < n$ and $\deg B < \deg f_0$.
Now
$$L(A(x))\circ L(f_0 (x_0)) + L(B(x))\circ L(\f_W (x_0) ) = x_0$$
$$L(A(x)) \circ L( \sum_{i=1}^{m-1} f_i (x_i)) = L (g)$$ for some $g=\sum_{i=1}^{m-1} g_i$ where $g_i\in k [x_i]$.
So $$L(A(x))\circ L(f) + L(B(x))\circ L(\f_W) = x_0+ L(g).$$

We have $$ L(A(x))\circ \ell (f) \equiv x_{00} +\ell(g)\pmod{\langle \bar{Q}\rangle}.$$
Note that  $\ell (g)\in S_{11}$.
By Lemma~\ref{L} there is some $f'\in S_1$ such that $L(A)\circ \ell(f) - f'\in \langle \bar{Q}\rangle$ and $f'\equiv_q 0 \pmod{\mathcal{G}}$.  So put $\ell_0=-\ell (g)$.
Then $f'\equiv x_{00} - \ell_0 \pmod{\langle \mathcal{G}\rangle}$, and since $f'$ and $x_{00} - \ell_0$ are both in $S_1$, we have
$f'=x_{00} - \ell_0$.
Let $\ell_1 \in S_{11}$ such that $\ell_0^q \equiv_q \ell_1 \pmod{\bar{Q}_{1}}$.
Then $x_{01}\equiv_q x_{00}^q\equiv_q \ell_0^q \equiv_q\ell_1
\pmod{\bar{\mathcal{G}}}$, and inductively we have $x_{0i}\equiv_q \ell_i$ for some linear form $\ell_i\in S_{11}$, with $\ell_{i-1}^q \equiv_q \ell_i \pmod{\bar{Q}_{1}}$.
$\Box$

When the condition in Lemma~\ref{xn-1} is satisfied, $x_{0i}\equiv_q \ell_i$ for some $\ell_i\in S_{11}$.   Substituting he variable $x_{0j}$ by $\ell_j$, for $j=0,\ldots,n'-1$, gives reduction from  $S_1\cap V_{\bar{\mathcal{G}},q}$ to $S_{11}\cap V_{\bar{\mathcal{G}},q}$.  More explicitly, for $g\in S_1$, write $g=g_0 + g_1$ where $g_0$ is a linear form in
$x_{00}$, ..., $x_{0\ n'-1}$, and $g_1\in S_{11}$.  Then
$g \equiv_1 g'$ where $g'=g_0 (\ell_0,\ldots,\ell_{n'-1}) + g_1\in S_{11}$.
Therefore for all $g\in S_1\cap V_{\bar{\mathcal{G}},q}$, there is some $g'\in S_{11}$ such that $0\equiv_q g \equiv_1 g' \pmod{\bar{\mathcal{G}}}$.
A similar condition will give reduction from  $S_{11}\cap V_{\bar{\mathcal{G}},q}$ to $S_{12}\cap V_{\bar{\mathcal{G}},q}$, and so on.  This leads to the following definition.

We say that $\mathcal{F}$ is {\em reducible} for $W$ if
for $i=0,\ldots, m-2$, either $V_{\bar{\mathcal{G}},q}\cap S_{1i}=V_{\bar{\mathcal{G}},q}\cap S_{1 i+1}$, or else there is a $k'$-linearized polynomial of the form $L(f_i)$ with $f_i=\sum_{j=i}^{m-1} g_{ij}$,  $g_{ij} \in k[x_j]$ of degree less than $n'$, for $j=i,\ldots,m-1$,
and $\ell(f_i)\in V_{\bar{\mathcal{G}},q}\cap S_{1i}$ and the GCD of $g_{ii}$ and $\f(x_i)$ is 1.

In particular if $\f_W$ is irreducible over $k'$ then the GCD of every nonzero polynomial of degree less than $n'=\deg\f_W$ is relatively prime to $\f_W$.  Therefore we have the following:

\begin{lemma}
If $\f_W$ is irreducible over $k'$ then $\mathcal{F}$ is reducible for $W$.
\end{lemma}

\begin{theorem}
\label{reducible-linearized}
Suppose $\mathcal{F}$ is a finite set of $k'$-linearized polynomials, and $W$ is a $\tau$-invariant subspace of $k$ where $\tau$ is the Frobenius map over $k'$. Let $\bar{\mathcal{G}}=\bar{\mathcal{F}}\cup\bar{Q}$.
If $\mathcal{F}$ is reducible for $W$, then $d_{\bar{\mathcal{G}}} \le (q-1)m+1$.  Moreover a basis of $Z_{W}(\mathcal{F})$ can be constructed in time
$(n' m)^{O(q)}$ where $n'=\deg\f_W$.
\end{theorem}

\begin{theorem}
\label{weil-linearized}
Suppose $\mathcal{F}$ is a finite set of $k'$-linearized polynomials of maximum degree $d=q^c$ for some $c > 0$.
Let $\mathcal{F}'$ be the Weil descent system of $\mathcal{F}$ with respect to a $k/k'$ basis, and $\mathcal{F}'_1=\mathcal{F}'\cup\{ x_{ij}^q-x_{ij}: i=0,\ldots,m-1, j=0,\ldots,n-1\}$.
If $\mathcal{F}$ is reducible for $k$, then $d_{\mathcal{F}'_1}\le \max( (q-1)m+1, qd)$.
\end{theorem}

\ \\{\bf Example}  Consider the case where $\mathcal{F}$ consists of a bivariate linearized polynomial
\begin{eqnarray*}
F(x,y)& = & ax^{q^2}+bx^{q} + c x + uy^{q^2} + v y^q + w y\\
  & = & L (ax^2 + bx + c) + L(uy^2+vy+w),
\end{eqnarray*}
with $a,b,c,u,v,w\in k=\F_{q^n}$.  By Lemma~\ref{xn-1} (with $f=ax^2+bx+c+uy^2+vy+w$), if either $GCD(ax^2 + bx + c,x^n-1) = 1$ or $GCD(uy^2+vy+w, y^n-1) = 1$, then $\mathcal{F}$ is reducible for $k$.  By Theorem~\ref{weil-linearized}, $d_{\mathcal{F}'_1}\le 2q$.  $\Box$

Since $\bar{\mathcal{G}}\subset V_{\mathcal{G}, d}$, Theorem~\ref{weil-linearized} follows from Theorem~\ref{weil} and Theorem~\ref{reducible-linearized}.  The rest of this section is devoted to the proof of Theorem~\ref{reducible-linearized}.

\subsection{Proof of Theorem~\ref{reducible-linearized}}
\begin{lemma}
\label{reducexij}
Suppose $\mathcal{F}$ is reducible for $W$. For $i=0,\ldots,m-2$, if $V_{\bar{\mathcal{G}},q}\cap S_{1i}\neq V_{\bar{\mathcal{G}},q}\cap S_{1 i+1}$, then $x_{ij} \equiv_q \gamma_{ij}\pmod{\bar{\mathcal{G}}}$ for some linear form $\gamma_{ij}\in S_{1\ m-1}$, for $j=0,\ldots,n'-1$;  moreover $\gamma_{ij}^q \equiv_q \gamma_{i\ j+1} \pmod{\bar{Q}_{m-1}}$ for $j=0,\ldots, n'-2$.
\end{lemma}

\ \\{\bf Proof}
For $i=0,\ldots,m-2$, if $V_{\bar{\mathcal{G}},q}\cap S_{1i}\neq V_{\bar{\mathcal{G}},q}\cap S_{1 i+1}$, then there is a $k'$-linearized polynomial of the form $L(f_i)$ with $f_i=\sum_{j=i}^{m-1} g_{ij}$,  where $g_{ij} \in k[x_j]$ of degree less than $n'$, for $j=i,\ldots,m-1$, $\ell(f_i)\in V_{\bar{\mathcal{G}},q}\cap S_{1i}$ and the GCD of $g_{ii}$ and $\f (x_i)$ is 1.

By Lemma~\ref{xn-1} we have the following:
for $j=0,\ldots,n'-1$, $x_{ij} \equiv_q \ell_{ij} \pmod{\bar{\mathcal{G}}}$ for some linear form $\ell_{ij} \in S_{1\ i+1}$, moreover $\ell_{ij}^q \equiv_q \ell_{i \ j+1}\pmod{\bar{Q}_{i+1}}$.  From this it is easy to see by induction (proceeding from $i= m-2$ to $i=0$) that
$x_{ij} \equiv_q \gamma_{ij}$ for some linear form $\gamma_{ij}\in S_{1\ m-1}$, moreover $\gamma_{ij}^q \equiv_q \gamma_{i\ j+1} \pmod{\bar{Q}_{m-1}}$ for $i=0,\ldots,m-1$, $j=0,\ldots, n'-2$. $\Box$

\begin{lemma}
\label{reduce-G-bar}
Let $\mathcal{N}=\{ i\in\{0,\ldots, m-2\}: V_{\bar{\mathcal{G}},q}\cap S_{1i}\neq V_{\bar{\mathcal{G}},q}\cap S_{1 i+1}\}$.
Let $\Gamma=\{ x_{ij}-\gamma_{ij}: \gamma_{ij}\in S_{1\ m-1}, x_{ij} \equiv_q \gamma_{ij}\pmod{\bar{\mathcal{G}}}, i\in\mathcal{N}, j=0,\ldots,n'-1\}$.
Let $H_{\bar{\mathcal{F}}} = \{ \ell(h_f): f\in\bar{\mathcal{F}}\}$.
Then there exist $H_1=  \{ \ell (h_{ij}) :i=0,\ldots,m-2, j=0,\ldots,n'-1\}$ where $h_{ij}\in k[x_{m-1}]$ with $\deg h_{ij} < n'$ such that letting $H= H_{\bar{\mathcal{F}}}\cup H_1$, then
$\Gamma\cup H\subset V_{\bar{\mathcal{G}},q}$,
$\bar{\mathcal{G}}=\bar{\mathcal{F}}\cup\bar{Q} \subset V_{H\cup\bar{Q}_{m-1}\cup\Gamma,q}$,
$\langle \bar{\mathcal{G}} \rangle = \langle H\cup \bar{Q}_{m-1}\cup \Gamma \rangle$.
\end{lemma}

\ \\{\bf Proof}
By Lemma~\ref{reducexij}, $\Gamma\subset V_{\bar{\mathcal{G}},q}$.
For
$f\in\bar{\mathcal{F}}\subset S_1$, let $f'\in S_{1\ m-1}$ be obtained from $f$ by substituting $x_{ij}$ with $\gamma_{ij}$ for $i=0,\ldots,m-2$, $j=0,\ldots,n'-1$. Then $f'\equiv_1 f \pmod{\Gamma}$, and $f'= \ell (h_f)$ for some $h_f\in k[x_{m-1}]$. We have
$\ell(h_f)\in V_{\Gamma,1}\subset V_{\bar{\mathcal{G}},q}$, and $f\in V_{H_{\bar{\mathcal{F}}}\cup \Gamma, 1}$ where
$H_{\bar{\mathcal{F}}} = \{ \ell(h_f): f\in\bar{\mathcal{F}}\}$.

For $i=0,\ldots, m-2$ and $j=0,\ldots, n'-2$, $x^q_{ij}-x_{i \ j+1}\in \bar{Q}$.   Lemma~\ref{reducexij} implies that
$x^q_{ij}-x_{i \ j+1}\equiv_q \gamma^q_{ij}-\gamma_{i \ j+1} \equiv_q \ell (h_{ij}) \pmod{\Gamma\cup \bar{Q}_{m-1} }$ for some $h_{ij}\in k[x_{m-1}]$.

For $i=0,\ldots, m-2$, $x_{i\ n'-1}^q - \ell (\g_W (x_i))\in \bar{Q}$.
Let $\ell (\g_W (x_i))=\sum_{j=0}^{n'-1} a_{ij} x_{ij}$ with $a_{ij}\in k$
Lemma~\ref{reducexij} implies that
$$x_{i \ n'-1}^q - \ell (\g_W (x_i)) \equiv_q \gamma_{i \ n'-2}^q - \sum_{j=0}^{n'-1} a_{ij} \gamma_{ij}\pmod{\Gamma }.$$
Since $\gamma_{i \ n'-2}^q \equiv_q \gamma_{i \ n'-1} \pmod{\bar{Q}_{m-1}}$, we have
$$ \gamma_{i \ n'-2}^q - \sum_{j=0}^{n'-1} a_{ij} \gamma_{ij}\equiv_q \ell(h_{i\ n'-1}) \pmod{\bar{Q}_{m-1}} $$ with
$h_{i\ n'-1}\in k[x_{m-1}]$ of degree less than $n'$.

To summarize, we have
$$x^q_{ij}-x_{i \ j+1}\equiv_q \ell (h_{ij}) \pmod{\Gamma\cup \bar{Q}_{m-1} }$$
for $i=0,\ldots, m-2$ and $j=0,\ldots, n'-2$,
and
$$x_{i \ n'-1}^q - \ell (\g_W (x_i)) \equiv_q \ell(h_{i\ n'-1})\pmod{\Gamma\cup \bar{Q}_{m-1} }$$
for $i=0,\ldots, m-2$.
Let $H_1=  \{ \ell (h_{ij}) :i=0,\ldots,m-2, j=0,\ldots,n'-1\}$.  It follows that
$\bar{Q}\subset V_{H_1\cup\Gamma\cup\bar{Q}_{m-1},q}$ and on the other hand
$H_1\subset V_{\bar{Q}\cup\Gamma, q}$, and since $\Gamma\subset V_{\bar{\mathcal{G}},q}$, we have
$H_1\subset V_{\bar{\mathcal{G}},q}$.

Let $H=H_{\bar{\mathcal{F}}}\cup H_1$.
Then we conclude that
$\Gamma\cup H\subset V_{\bar{\mathcal{G}},q}$, and on the other hand
$\bar{\mathcal{G}}=\bar{\mathcal{F}}\cup\bar{Q} \subset V_{H\cup\bar{Q}_{m-1}\cup\Gamma,q}$.  In particular, we have
$\langle \bar{\mathcal{G}} \rangle = \langle H\cup \bar{Q}_{m-1}\cup \Gamma \rangle$.  $\Box$

Note that $H\cup \bar{Q}_{m-1}\subset k [\hat{x}_{m-1}]$ where $\hat{x}_{m-1} = x_{m-1\ 0},\ldots, x_{m-1 \ n'-1}$.

\begin{lemma}
\label{reduce-to-g}
Let $H$ be as in Lemma~\ref{reduce-G-bar}.
Suppose $H=\{\ell(h_i): i=1,\ldots,s\}$ and let $h_0 = \f_W (x_{m-1})$.
Let $g$ be the GCD of $h_i$, $i=0,\ldots,s$.
Then $\langle \bar{\mathcal{G}} \rangle= \langle \Gamma\cup \{\ell(g)\} \cup\bar{Q}_{m-1} \rangle$, moreover $\Gamma\cup \{\ell(g)\} \cup\bar{Q}_{m-1} \subset V_{\bar{\mathcal{G}},q}$.
\end{lemma}

\ \\{\bf Proof}
We have $g=\sum_{i=0}^s a_i h_i$ with $a_i\in k [ x_{m-1}]$, so
$$ L(g)= \sum_i L(a_i)\circ L(h_i).$$  So
$$\ell(g)\equiv \sum_i L(a_i) \circ \ell (h_i) \mod {\bar{Q}_{m-1}}.$$
Apply Lemma~\ref{L} to $H\cup \bar{Q}_{m-1}\subset k [\hat{x}_{m-1}]$ it follows that there is $h'_i\in S_{1\ m-1}$ such that $h'_i\equiv_q 0 \pmod{H\cup\bar{Q}_{m-1}}$ and
$L(a_i) \circ \ell (h_i) \equiv h'_i \pmod{\langle \bar{Q}_{m-1} \rangle}$.  So
$$\ell(g) \equiv \sum_i h'_i \equiv_q 0 \pmod{H\cup\bar{Q}_{m-1}}.$$  Since $\ell(g)$ and $h'_i$ are all in $S_{1\ m-1}$, we have
$$\ell(g) =\sum_i h'_i \equiv_q 0 \pmod{H\cup\bar{Q}_{m-1}},$$
in particular, $\ell(g)\equiv_q 0 \pmod{\bar{\mathcal{G}}}$.
It follows that $$L(\f_W)=L(h_0)\in\langle \{ L(h_i): i=0,\ldots,s\}\rangle = \langle L(g) \rangle.$$
Under the isomorphism from $k[\hat{x}_{m-1}]/\langle \bar{Q}_{m-1}\rangle \to k[ x_{m-1}]/\langle L (\f_W (x_{m-1}))\rangle$, $\ell(h_i)$ corresponds to $L(h_i)$, hence the ideal generated by
$H\cup\bar{Q}_{m-1}$ corresponds to the ideal generated by $L(g)$.

Since by Lemma~\ref{reduce-G-bar}, $\langle \bar{\mathcal{G}} \rangle= \langle \Gamma\cup H\cup\bar{Q}_{m-1} \rangle$, it follows that
$\langle \bar{\mathcal{G}} \rangle= \langle \Gamma\cup \{\ell(g)\} \cup\bar{Q}_{m-1} \rangle$.  Moreover from the discussion above we have $\Gamma\cup \{\ell(g)\} \cup\bar{Q}_{m-1} \subset V_{\bar{\mathcal{G}},q}$.  $\Box$

Under the isomorphism from $k[\hat{x}_0,\ldots,\hat{x}_{m-1}]/\langle \bar{Q}\rangle \to k[x_0,\ldots,x_{m-1}]/\langle Q\rangle$, $x_{i0}-\gamma_{i0}$ corresponds to $x_i - L(g_i)$ where $\ell (g_i)=\gamma_{i0}$ for $i\in\mathcal{N}$.
Under the isomorphism the ideal determined by $\bar{\mathcal{G}}$ corresponds to the ideal determined by $\mathcal{F}\cup Q $.
Since, by Lemma~\ref{reduce-to-g}, $\langle \bar{\mathcal{G}} \rangle= \langle \Gamma\cup \{\ell(g)\} \cup\bar{Q}_{m-1} \rangle$ and $g | \f_W$,
it follows
that $\langle \mathcal{F}\cup Q \rangle$ is generated by $L(g)$ and $x_i - L(g_i)$ where $i\in\mathcal{N}$.
By Lemma~\ref{reduce-to-g}  $\Gamma\cup \{\ell(g)\} \cup\bar{Q}_{m-1} \subset V_{\bar{\mathcal{G}},q}$, it follows from Proposition 2.3 of \cite{HKYY} that $\ell(g)$ and $\gamma_{i0}$, hence $L(g)$ and $x_i-L(g_i)$ can be constructed in time $(mn')^{O(q)}$ time.
From this a basis of
$Z_W (\mathcal{F})$ over $k'$ can be easily written down.

It is easy to see that if $f\in k[\hat{x}_{m-1}]$ and $f\in \langle \{\ell(g)\} \cup\bar{Q}_{m-1} \rangle$ then
$f\equiv_{\deg f +1} \ell(g) f_1 \pmod{\bar{Q}_{m-1}}$ for some $f_1\in k [\hat{x}_{m-1}]$.
Suppose $f\in \langle \bar{\mathcal{G}} \rangle$.  Then $f\equiv_{\deg f} f_1\pmod{\bar{Q}}$ where the degree of $x_{ij}$ in $f_1$ is less than $q$ for all $i,j$.  Using $x_{ij}\equiv \gamma_{ij}\pmod{\Gamma}$, we have $f_1 \equiv_{\deg f_1} h \pmod{\Gamma\cup\bar{Q}_{m-1}}$ where $h\in k[\hat{x}_{m-1}]$.  It follows that $h\in\langle \{\ell(g)\}\cup \bar{Q}_{m-1}\rangle$, hence $h\equiv_{\deg h +1} \ell(g)h_1\pmod{\bar{Q}_{m-1}}$, so
$h\equiv_{\deg h +1} 0 \pmod{\{\ell(g)\}\cup \bar{Q}_{m-1}}$.
If $\deg f > (q-1) m$, then $\deg f > \deg f_1$,
and since $ \Gamma\cup \{\ell(g)\} \cup\bar{Q}_{m-1} \subset V_{\bar{\mathcal{G}},q}$, we conclude that
$f\in V_{\bar{\mathcal{G}},\deg f}$. Therefore $d_{\bar{\mathcal{G}}} \le (q-1)m+1$.  Theorem~\ref{reducible-linearized} follows.

\end{document}